\documentclass{amsart}
\makeatletter

\def\@setabstract{\@setabstracta \let\@setabstract\relax}

\renewcommand{\author}[2][]{%
  \ifx\@empty\authors
    \gdef\authors{#2}%
  \else
    \g@addto@macro\authors{\and#2}%
    \g@addto@macro\addresses{\author{}}%
  \fi
  \@ifnotempty{#1}{%
    \ifx\@empty\shortauthors
      \gdef\shortauthors{#1}%
    \else
      \g@addto@macro\shortauthors{\and#1}%
    \fi
  }%
}
\edef\author{\@nx\@dblarg
  \@xp\@nx\csname\string\author\endcsname}
\renewcommand\@settitle{\noindent{\huge\bfseries\@title\vskip1.5em}}
\renewcommand\@setauthors{\noindent{\Large \andify\authors \authors}}

\renewcommand{\maketitle}{
  \par
  \@topnum\z@ 
  \thispagestyle{empty}
  \uppercasenonmath\shorttitle
  \ifx\@empty\shortauthors \let\shortauthors\shorttitle
  \else \andify\shortauthors \fi
  \begingroup
  \@maketitle
  \toks@\@xp{\shortauthors}\@temptokena\@xp{\shorttitle}%
  \edef\@tempa{%
    \@nx\markboth{\@nx\MakeUppercase{\the\toks@}}{\the\@temptokena}}\@tempa
  \endgroup
  \c@footnote\z@
  \def\do##1{\let##1\relax}%
  \do\maketitle \do\@maketitle \do\title \do\@xtitle \do\@title
  \do\author \do\@xauthor \do\address \do\@xaddress
  \do\email \do\@xemail \do\curraddr \do\@xcurraddr
  \do\dedicatory \do\@dedicatory \do\thanks \do\thankses
  \do\keywords \do\@keywords \do\subjclass \do\@subjclass
}
\renewenvironment{abstract}{
  \ifx\maketitle\relax
    \ClassWarning{\@classname}{Abstract should precede
      \protect\maketitle\space in AMS documentclasses; reported}
  \fi
  \global\setbox\abstractbox=\vtop \bgroup
    \normalfont
    \list{}{\labelwidth\z@
      \leftmargin3pc \rightmargin\leftmargin
      \listparindent\normalparindent \itemindent\z@
      \parsep\z@ \@plus\p@
      
    }
    \item[\hskip\labelsep\scshape\abstractname.]
}{
  \endlist\egroup
  \ifx\@setabstract\relax \@setabstracta \fi
}

\def\@setaddresses{\par
  \nobreak \begingroup
\footnotesize
\setlength{\parindent}{0pt} %
  \def\author##1{\nobreak\addvspace\bigskipamount}%
  \interlinepenalty\@M
  \def\address##1##2{\begingroup
    \par\addvspace\bigskipamount\indent
    \@ifnotempty{##1}{(\ignorespaces##1\unskip) }%
    {\scshape\ignorespaces##2}\par\endgroup}%
  \def\curraddr##1##2{\begingroup
    \@ifnotempty{##2}{\nobreak\indent{\itshape Current address}%
      \@ifnotempty{##1}{, \ignorespaces##1\unskip}\/:\space
      ##2\par}\endgroup}%
  \def\email##1##2{\begingroup
    \@ifnotempty{##2}{\nobreak\indent{\itshape E-mail address}%
      \@ifnotempty{##1}{, \ignorespaces##1\unskip}\/:\space
      \ttfamily##2\par}\endgroup}%
  \def\urladdr##1##2{\begingroup
    \@ifnotempty{##2}{\nobreak\indent{\itshape URL}%
      \@ifnotempty{##1}{, \ignorespaces##1\unskip}\/:\space
      \ttfamily##2\par}\endgroup}%
  \addresses
  \endgroup
}

\def\@secnumfont{\bfseries}
\renewcommand\section{\@startsection {section}{1}{\z@}%
  {-3.5ex \@plus -1ex \@minus -.2ex}{2.3ex \@plus.2ex}{\normalfont\Large\bfseries}}
\renewcommand\subsection{\@startsection{subsection}{2}{\z@}%
  {-3.25ex\@plus -1ex \@minus -.2ex}{1.5ex \@plus .2ex}{\normalfont\large\bfseries}}
\renewcommand\subsubsection{\@startsection{subsubsection}{3}{\z@}%
  \z@{-.5em}{\normalfont\normalsize\bfseries}}
\makeatother

\newtheoremstyle{theoremlike}{3pt}{3pt}{\itshape}{}{\bfseries}{.}{.5em}{}
\theoremstyle{theoremlike}
\newtheorem{theorem}{Theorem}[section]

\newtheorem{lemma}[theorem]{Lemma}

\newtheorem{definition}[theorem]{Definition}

\theoremstyle{definition}
\theoremstyle{remark}

\usepackage{amsmath}
\usepackage{amsthm}
\usepackage{amssymb}
\usepackage{color}
\usepackage{dsfont}
\usepackage{graphicx}
\usepackage{graphics}

\usepackage[margin=1.in]{geometry}

\usepackage{caption}
\usepackage{subcaption}
\begin{document}

\title[On Third-Order Limiter Functions for Finite Volume Methods]{On Third-Order Limiter Functions for Finite Volume Methods}

\author{Birte \textsc{Schmidtmann}$^{a) ,}$*}
\address{$^{a)}$ Center for Computational Engineering Science \\ RWTH Aachen University \\
D-52062 Aachen \\ Germany}
\email{schmidtmann@mathcces.rwth-aachen.de}

\author{R\'emi \textsc{Abgrall}$^{b)}$}
\address{$^{b)}$ Institut f\"ur Mathematik $\&$ Computational Science\\ Universit\"at Z\"urich \\ CH-8057 Z\"urich \\ Switzerland}
\email{remi.abgrall@math.uzh.ch}
\author{Manuel \textsc{Torrilhon}$^{a)}$}

\begin{abstract}
	In this article, we propose a finite volume limiter function for a reconstruction on the three-point stencil. Compared to classical limiter functions in 
	the MUSCL framework, which yield $2^\text{nd}$-order accuracy, the new limiter is $3^\text{rd}$-order accurate for smooth solution. 
	In an earlier work, such a $3^\text{rd}$-order limiter function was proposed and showed successful results \cite{CT2009}. However,
	it came with unspecified parameters. We close this gap by giving information on these parameters.
\end{abstract}

\begin{center}
\maketitle
\end{center}
%
\section{Introduction}\label{sec:introduction}
	%
	We consider the numerical approximation of hyperbolic conservation laws of the form
	\begin{subequations}	
	\label{eq:consLaw}
		\begin{align}
			u_t + (f(u))_x &= 0, \\
			u(x,0) &= u_0 (x),
		\end{align}
	\end{subequations}	
	where $u=(u_1,\ldots,u_s)^T$ and the Jacobian matrix $A(u) = \partial f / \partial u$ has $s$ real eigenvalues. 
	In this work, we restrict our discussion to the scalar 1D case $s=1$. We further assume $u_0(x)$ to be either periodic or to have compact support.\\
	On a regular computational grid with space intervals of size $\Delta x$, let $x_i$ denote the position of the cell centers. The control cells
	are defined by $C_i = [x_{i-1/2}, x_{i+1/2}]$, where $x_{i\pm 1} = x_i \pm \Delta x$.\\
	The solution of Eq. \eqref{eq:consLaw} is approximated by the cell averages $\bar{u}_i^n = \frac{1}{\Delta x} \int_{C_i} u(x,t^n) dx$ which are updated with 
	the finite volume (FV) formulation of Eq. \eqref{eq:consLaw} given by
	\begin{align}
		\label{eq:approxConsLaw}
		\frac{d \bar u_i}{d t} = -\frac{1}{\Delta x} \left( \hat f_{i+1/2} - \hat f_{i-1/2}\right).	
	\end{align}
	The numerical flux function $\hat f_{i\pm 1/2} = f(u(x_{i \pm 1/2},t))$ results from Eq. \eqref{eq:consLaw} by integrating over $C_i$. The aim is to define an 
	update rule for the new time step $t^{n+1} = t^n + \Delta t$ such that Eq. \eqref{eq:consLaw} is approximated with high order of accuracy. The main 
	challenge is to avoid the development of spurious oscillations near shocks and at the same time maintain high order accuracy at smooth extrema.
	\\	
	\\
	We are interested in a numerical scheme with the most compact stencil, using only information of the cell $C_i$ and its most direct neighbors $C_{i-1}$ and $C_{i+1}$. Classical 
	approaches based on this three-point-stencil, such as the MUSCL scheme, yield $2^{\text{nd}}$ order schemes \cite{VL1979, LeVeque2002}, however, we will present an update rule that yields 
	$3^{\text{rd}}$ order accuracy for smooth solutions.
	\\
	The key point is the definition of the numerical flux function $\hat f$ which depends on the left and right limiting values $u^{(\pm)}(x_{i\pm 1/2})$ at the 
	cell boundaries $x_{i\pm 1/2}$, cf. Fig \ref{fig:reconstruction}. These values are a priori not known and have to be reconstructed from the cell mean values $\bar{u}_i^n$. The focus of this 
	work is on the reconstruction procedure.
	%
\section{Theory}\label{sec:theory}
\subsection{Two Parameter Setting}\label{subsec:theory1}
\begin{figure}
		\begin{minipage}{0.49\textwidth}		
			\centering
			\includegraphics[width=\textwidth]{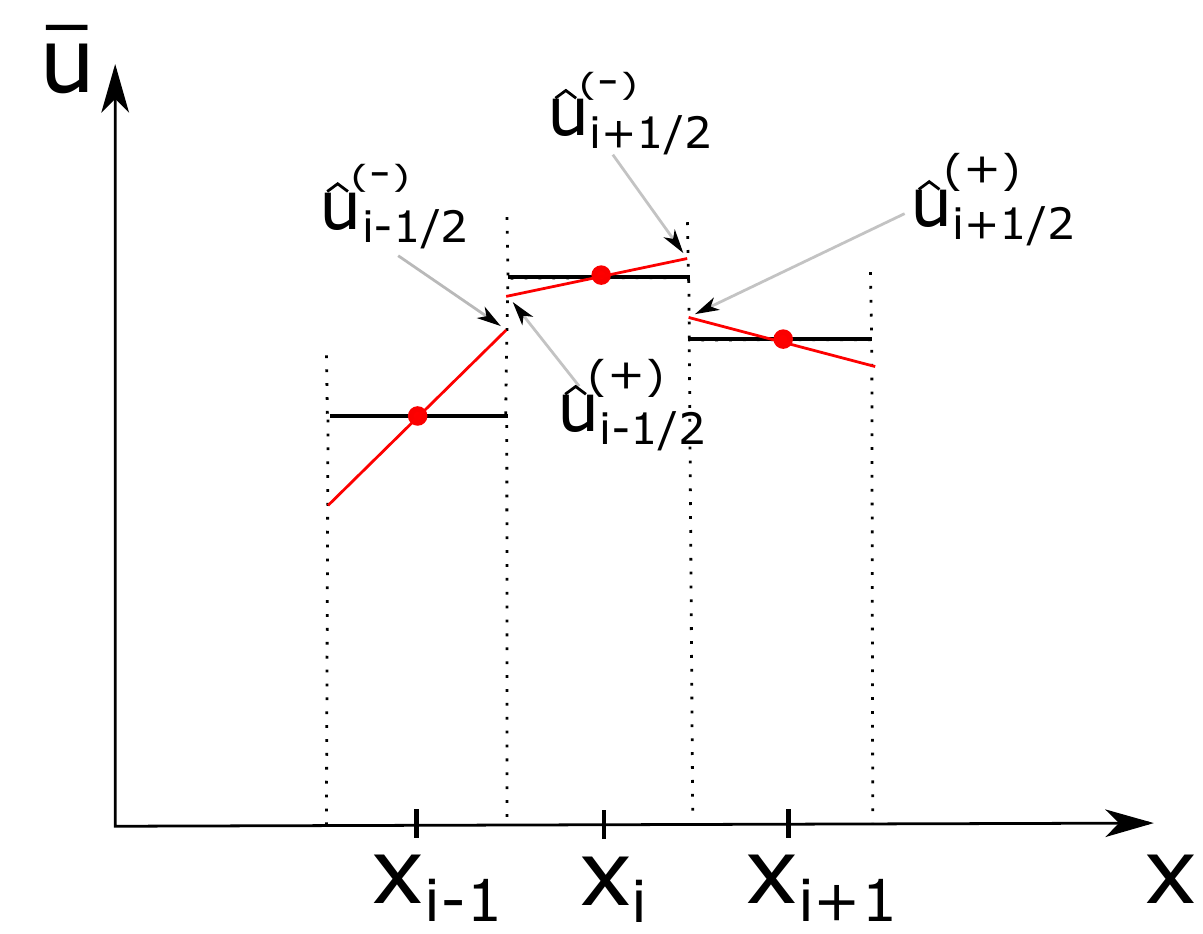}
		\end{minipage}
		\begin{minipage}{0.49\textwidth}
			\centering
			\includegraphics[width=\textwidth]{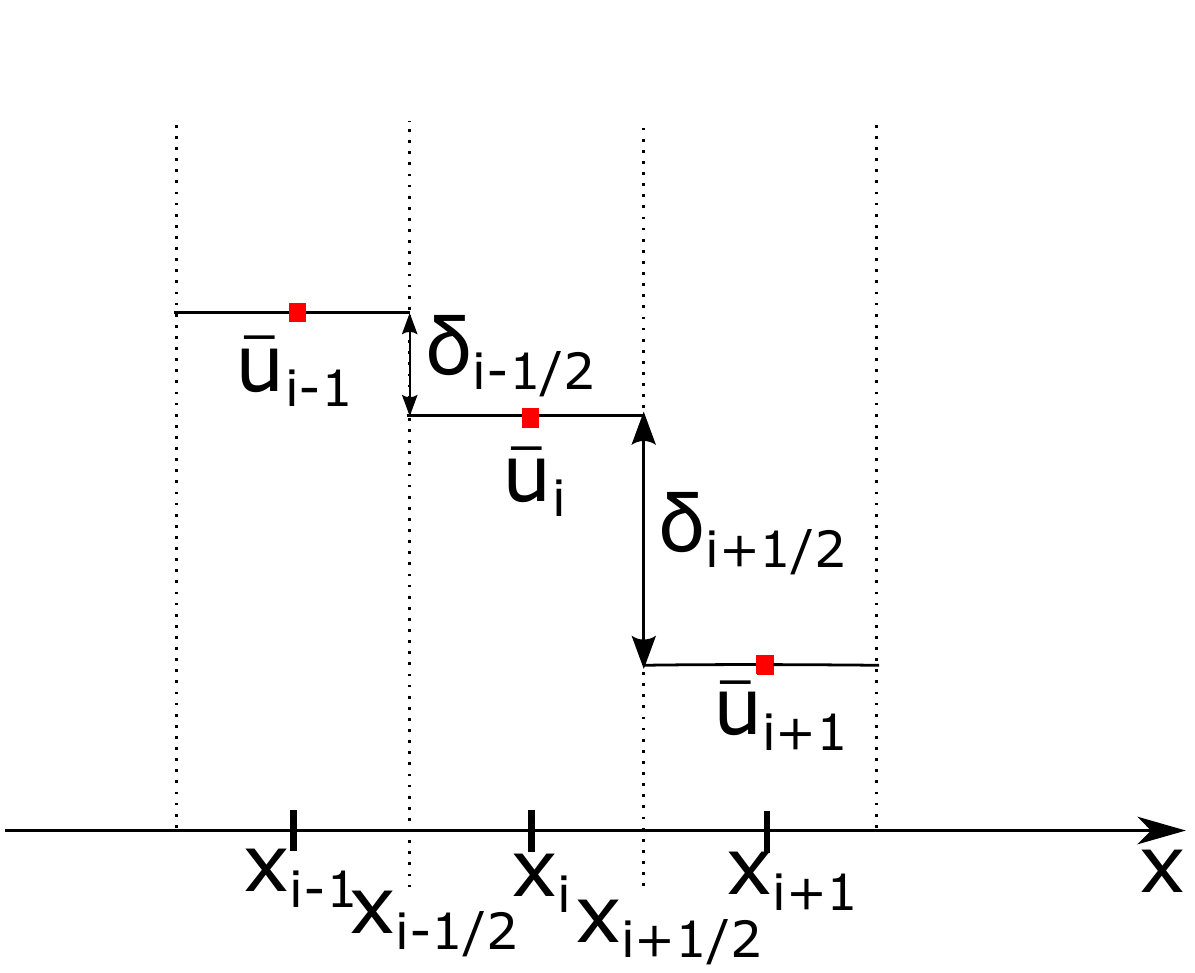}
		\end{minipage}
		\caption{Basic setting for the reconstruction of the interface values $u^{(\pm)}(x_{i\pm 1/2})$ on a 3-point-stencil.}
		\label{fig:reconstruction}
	\end{figure}

\begin{figure}
		\centering
		\includegraphics[scale=0.4]{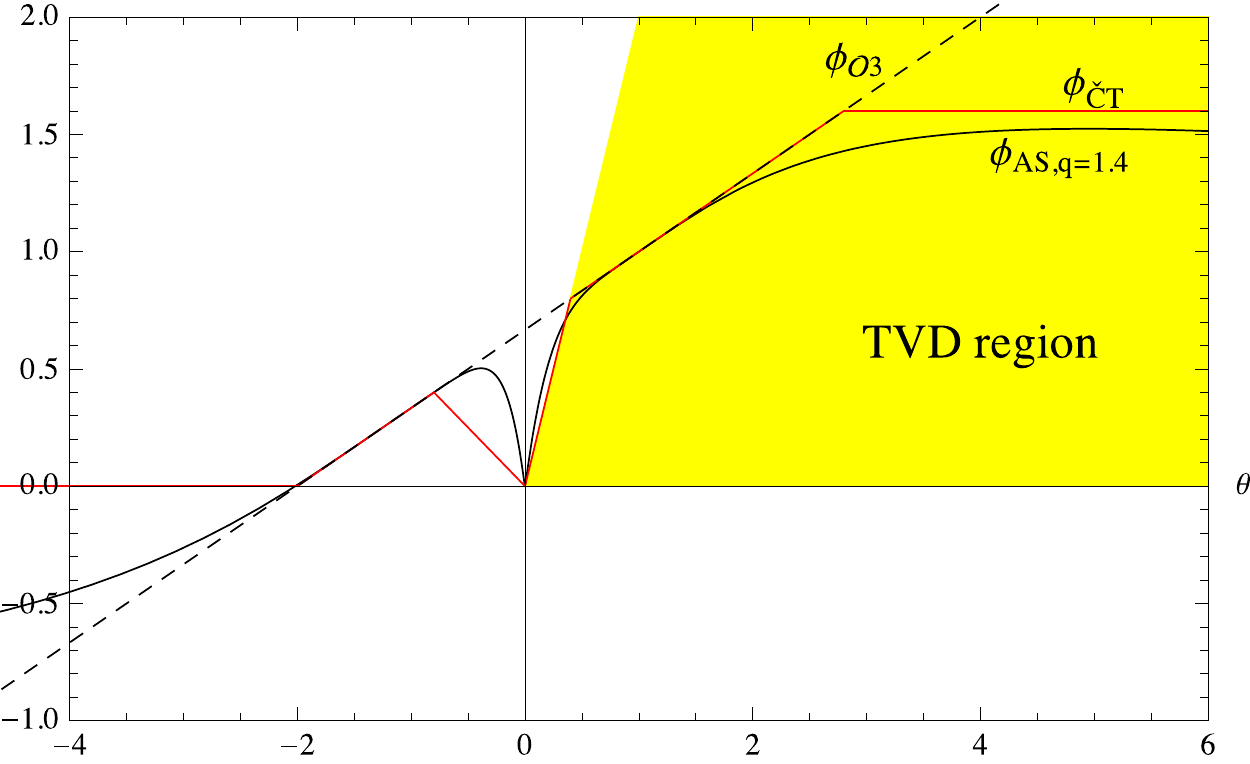}
		\caption{Alternative limiter matching the logarithmic limiter with $q=1.4$.}
		\label{fig:ASandCTlimiter}
	\end{figure}

\begin{figure}[ht!]
   	\centering
   		\begin{subfigure}{0.49\textwidth}
   			\centering
   			\includegraphics[width=\textwidth]{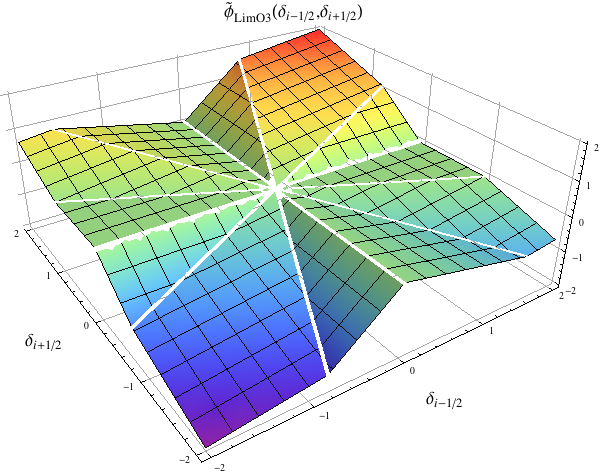}   		
	   		\caption{Alternative limiter function $\tilde \phi_{\text{Lim}\mathcal{O}3}$.}
			\label{fig:CadaLimiter2dFramework}
   		\end{subfigure}
   		\hfill
   		\begin{subfigure}{0.49\textwidth}
   			\centering
   			\includegraphics[width=\textwidth]{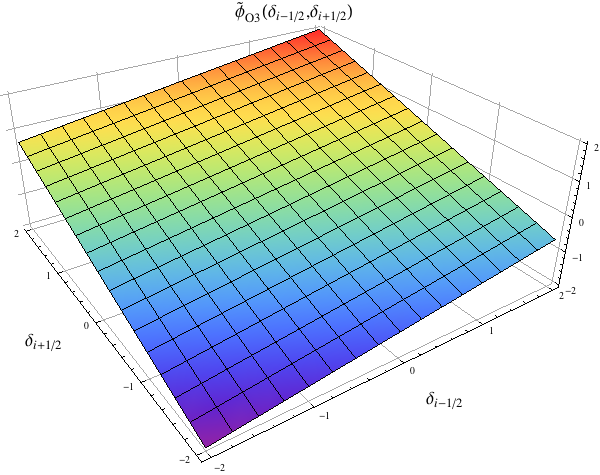}   		
	   		\caption{Full-third-order reconstruction $\tilde\phi_{\mathcal{O}3}$.}
			\label{fig:3rdOrderLimiter2dFramework}
   		\end{subfigure}
   		\caption{Different reconstruction functions in the two-parameter-framework.}
   \end{figure}

	%
	Considering the compact stencil $\{x_{i-1}, x_i, x_{i+1}\}$, we want to reconstruct the interface values at the cell boundaries $x_{i\pm 1/2}$ as 
	shown in Fig. \ref{fig:reconstruction}. For the cell $x_i$, we use the left and right interface values defined by
	\begin{subequations}
    	\label{eq:reconstructionNonStandard}
		\begin{align}
	 	  	u^{(-)}_{i+1/2} &= \bar u_i + \frac{1}{2}\; \phi(\theta_i)\delta_{i+1/2} \\
  		  	u^{(+)}_{i-1/2} &= \bar u_i - \frac{1}{2}\; \phi(\theta^{-1}_i)\delta_{i-1/2}.
		\end{align}
 	\end{subequations}
	%
	Here, $\phi$ is a non-linear limiter function depending on the local smoothness measure
	\begin{align*}
		\theta_i = \frac{\delta_{i-1/2}}{\delta_{i+1/2}},\quad \delta_{i+1/2}\neq 0 		
	\end{align*}
	with $\delta_{i+1/2}=\bar u_{i+1}-\bar u_i,\; \delta_{i-1/2} =\bar u_i-\bar u_{i-1}$, cf. Fig. \ref{fig:reconstruction}. In Eq. \eqref{eq:reconstructionNonStandard}, 
	the choice of $\phi(\theta_i)$ determines the order of accuracy of the reconstruction and therefore of the scheme.\\ \\
	\noindent There is a variety of schemes on the three-point stencil that obtain $2^{\text{nd}}$-order accuracy. These are the classical schemes, 
	which use the information of the three cells to compute a linear reconstruction function, see e.g. \cite{VL1979}. Indeed, the second-order reconstruction 
	$u^{(-)}_{i+1/2}=\bar u_i + \frac{\Delta x}{2}	\left(\frac{\bar u_{i+1}-\bar u_{i-1}}{2\Delta x} \right)$ can be rewritten
	in form of Eq. \eqref{eq:reconstructionNonStandard} with the limiter function $\phi (\theta) = \frac{1+\theta}{2}$.
	This limiter function has the property that $\phi(\theta^{-1})=\theta^{-1}\phi(\theta)$ holds and therefore, Eq. 
	\eqref{eq:reconstructionNonStandard} can be reduced to the standard formulation	
	\begin{subequations}
    	\label{eq:reconstruction}
		\begin{align}
	 	  	u^{(-)}_{i+1/2} &= \bar u_i + \frac{\Delta x}{2}\; \sigma_i \\
  		  	u^{(+)}_{i-1/2} &= \bar u_i - \frac{\Delta x}{2}\; \sigma_i,
		\end{align}
  	\end{subequations}
	with the downwind slope $\sigma_i = \phi (\theta)\delta_{i+1/2}$ (see e.g. \cite{LeVeque2002}). 
	The aim of this work is to introduce schemes which use the three-point stencil to achieve $3^{\text{rd}}$ order accurate
	reconstructions of the cell-interface values. One possibility is to construct a quadratic polynomial $p_i(x)$ in each cell. 
	Applying the computed polynomial to $x_{i\pm 1/2}$ yields the interface values 
	\begin{subequations}
		\begin{align}
			u^{(+)}_{i-1/2} = p_i(x_{i-1/2})\\
			u^{(-)}_{i+1/2} = p_i(x_{i+1/2}).
		\end{align}
	\end{subequations}
	Rewriting the interface values 	in the form \eqref{eq:reconstructionNonStandard} yields 
	\begin{align}
	\label{eq:3rdOrder}
	  \phi_{\mathcal{O}3}(\theta_i) = \frac{2+\theta_i}{3}.
	\end{align}
	This formulation results in a full third order scheme for smooth solutions, however, causes oscillations near shocks
	and discontinuities. Since this should be avoided, we introduce a limiter function $\tilde\phi$, that 
	applies the full $3^{\text{rd}}$ order reconstruction Eq. \eqref{eq:3rdOrder} at smooth parts of the solution and switches 
	to a lower order reconstruction formulation close to large gradients, shocks and discontinuities.
	\newline
	\newline
	The limiting function we will dwell upon in this paper is based on the local double logarithmic reconstruction function 
	of Artebrant and Schroll \cite{AS2005}. They present a limiter function $\phi_{AS}(\theta_i, q)$ which contains 
	an additional parameter $q <1$. This parameter significantly changes the reconstruction function. The authors state that $q=1.4$ is the best choice and for $q\to 0$, the logarithmic
	limiter function reduces to $\phi_{\mathcal{O}3}(\theta_i)$, Eq. \eqref{eq:3rdOrder}.
	\begin{align*}
		&\phi_{AS}(\theta_i, q) =\frac{2 p[(p^2-2 p\,\theta+1)log(p)-(1-\theta)(p^2-1)]}{(p^2-1)(p-1)^2},\\
		&p=p(\theta_i, q)=2\frac{|\theta_i|^q}{1+|\theta_i|^{2q}}.
	\end{align*}
	The drawback of $\phi_{AS}(\theta_i, q)$ is its complexity which makes the evaluation in each cell expensive and possibly instable.\\
	In \cite{CT2009}, {\v{C}}ada and Torrilhon develop a limiter function $\phi_{\text{Lim}\mathcal{O}3}(\theta_i)$ that resembles the properties 
	of $\phi_{AS}$ and reduces the computational cost. The alternative limiter function reads
	\begin{align*}
		\phi_{\text{Lim}\mathcal{O}3}(\theta_i)=\max\left(0,\min\left( \phi_{\mathcal{O}3}(\theta_i), \max\left(-\frac{1}{2}\theta_i,\min\left(2 \,\theta_i,\phi_{\mathcal{O}3}(\theta_i),1.6 \right)\right)\right)\right)
	\end{align*}
	and is shown in Fig. \ref{fig:ASandCTlimiter} together with $\phi_{AS}(\theta_i, 1.4)$ and $\phi_{\mathcal{O}3}(\theta_i)$.\\
	\\
	All reconstruction functions presented so far have non-zero values for $\theta < 0$, which means that they break
	with the total variation diminishing (TVD) property. The idea of keeping the non-zero part in the construction 
	of $\phi_{\text{Lim}\mathcal{O}3}(\theta)$ for $\theta \in [-2, 0]$ was to avoid the clipping of smooth extrema. Extrema clipping is the effect that
	occurs close to minima and maxima, where the normalized slopes $\delta_{i\pm 1/2}$ are of the same order of magnitude but have opposite signs, i.e. $\theta 
	\approx -1$. In this case, classical limiter functions that fully lie in the TVD region yield zero and thus $1^{\text{st}}$ order accuracy. 
	This effect is avoided including the non-zero part in $\phi_{\text{Lim}\mathcal{O}3}$. \\
	Another clipping phenomenon arises, if the discretization of a smooth function contains a zero slope, $\delta_{i-1/2}\approx 0$
	or $\delta_{i+1/2}\approx 0$. This leads to $\theta \to 0$ or $\theta\to \pm \infty$ and the interface values are approximated
	by the cell mean values, which yields a $1^\text{st}$ order scheme. 	This case shows, that we need a criterion that can differentiate between 
	smooth extrema and discontinuities. 
	We require this decision criterion to depend only on information available on the compact three-point stencil. Furthermore, 
	it has to detect cases when to switch to the $3^\text{rd}$ order reconstruction, Eq. \eqref{eq:3rdOrder}, in case of smooth extrema, even though one of 
	the normalized slopes is zero. This is the case if the non-zero slope is 'small', compared to the case of 
	a discontinuity. The main focus of this work is to determine what 'small' means and to define a switch function $\eta$.\\
	\newline
	From the discussion above, it is clear that $\eta$ has to explicitly depend on both normalized slopes $\delta_{i\pm1/2}$. The classical approach
	of considering the ratio $\theta_i$ of neighboring slopes is overly restrictive because part of the information is given away.
	This is why we reformulate the limiter functions $\phi$ in a two-parameter-framework and obtain the new formulation for the reconstructed interface 
	values
	\begin{subequations}
    	\label{eq:reconstructionNew}
		\begin{align}
	 	  	u^{(-)}_{i+1/2} &= \bar u_i + \frac{1}{2}\; \tilde \phi(\delta_{i-1/2},\delta_{i+1/2}), \\
  		  	u^{(+)}_{i-1/2} &= \bar u_i - \frac{1}{2}\; \tilde \phi(\delta_{i+1/2},\delta_{i-1/2}),
		\end{align}
  	\end{subequations}
   	where the limiter function in the two-parameter framework is defined by
 	\begin{align}
 		\label{eq:2paramLimiter}
 		\tilde{\phi}(\delta_{i-1/2}, \delta_{i+1/2}) = \phi(\delta_{i-1/2}/\delta_{i+1/2})\delta_{i+1/2}.
	\end{align}
	This formulation avoids the division by the normalized slope which can be close to zero and thus cause instabilities.\\
    In this setting, the full-third-order reconstruction, Eq. \eqref{eq:3rdOrder}, reads
	\begin{align}
	\label{eq:3rdOrder2param}
	  \tilde\phi_{\mathcal{O}3}(\delta_{i-1/2}, \delta_{i+1/2}) = \frac{2 \delta_{i+1/2} + \delta_{i-1/2}}{3}.
	\end{align}
   	Fig. \ref{fig:CadaLimiter2dFramework} shows the alternative limiter function $\tilde \phi_{\text{Lim}\mathcal{O}3}$ and the full-third-order reconstruction
   	$\tilde\phi_{\mathcal{O}3}$ in the two-parameter setting. 
	\\
    On the coordinate axis, where $\delta_{i-1/2}=0$, i.e. $\theta_i=0$, the limiter function $\tilde \phi_{\text{Lim}\mathcal{O}3}$ returns zero, 
    meaning that it yields a $1^{\text{st}}$ order method. The same holds for the coordinate axis where $\delta_{i+1/2}=0$, see Eq. \eqref{eq:2paramLimiter}.
 	For two consecutive slopes of approximately the same order of magnitude, i.e. around the diagonals, the $3^{\text{rd}}$ order reconstruction Eq. 
 	\eqref{eq:3rdOrder2param} is gained.	\\
 	\begin{figure}[ht!]
   	\centering
   		\begin{subfigure}{0.49\textwidth}
   			\centering
   			\includegraphics[width=0.9\textwidth]{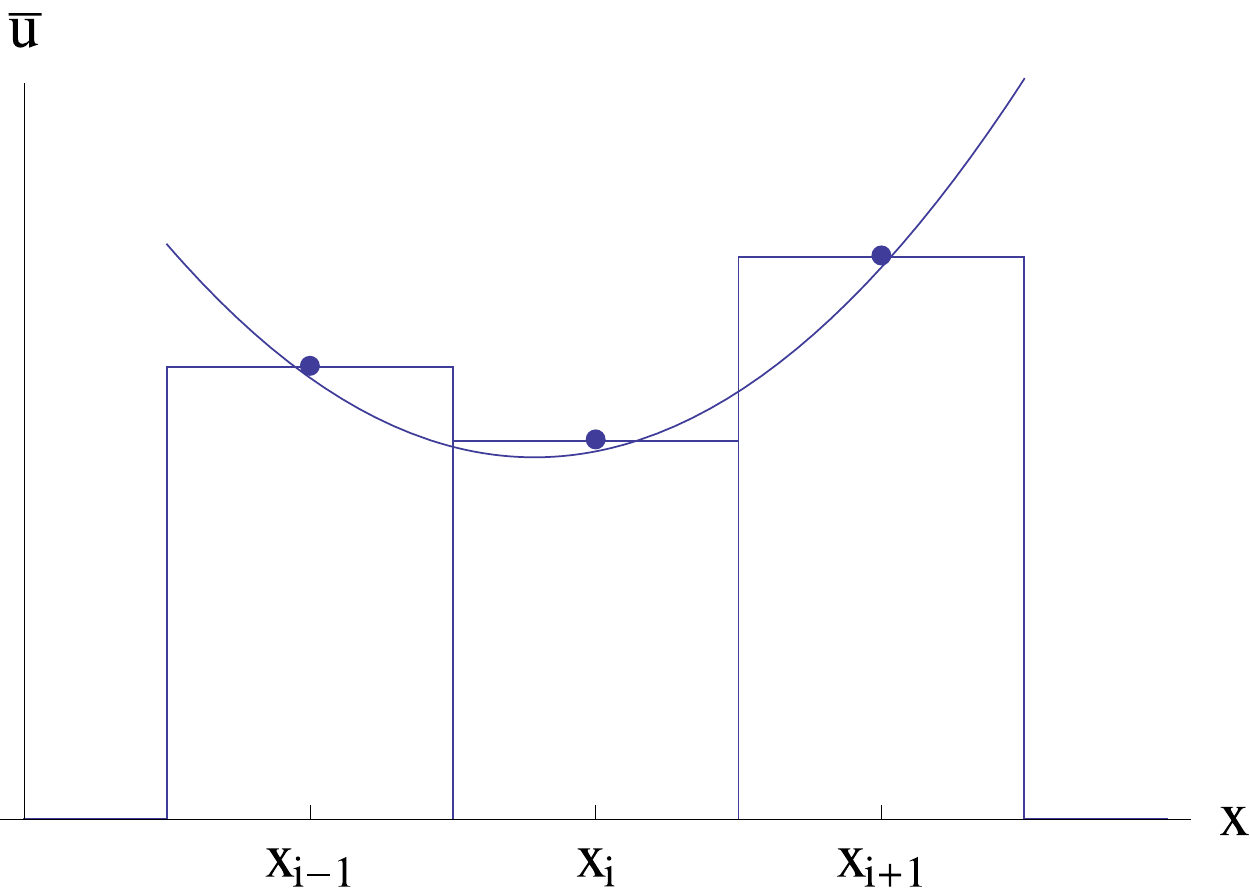}   		
	   		\caption{This situation is treated as a possible discontinuity: $\tilde \phi_{\text{Lim}\mathcal{O}3}\neq \tilde \phi_{\mathcal{O}3}$.}
   		\end{subfigure}
   		\hfill
   		\begin{subfigure}{0.49\textwidth}
   			\centering
   			\includegraphics[width=0.9\textwidth]{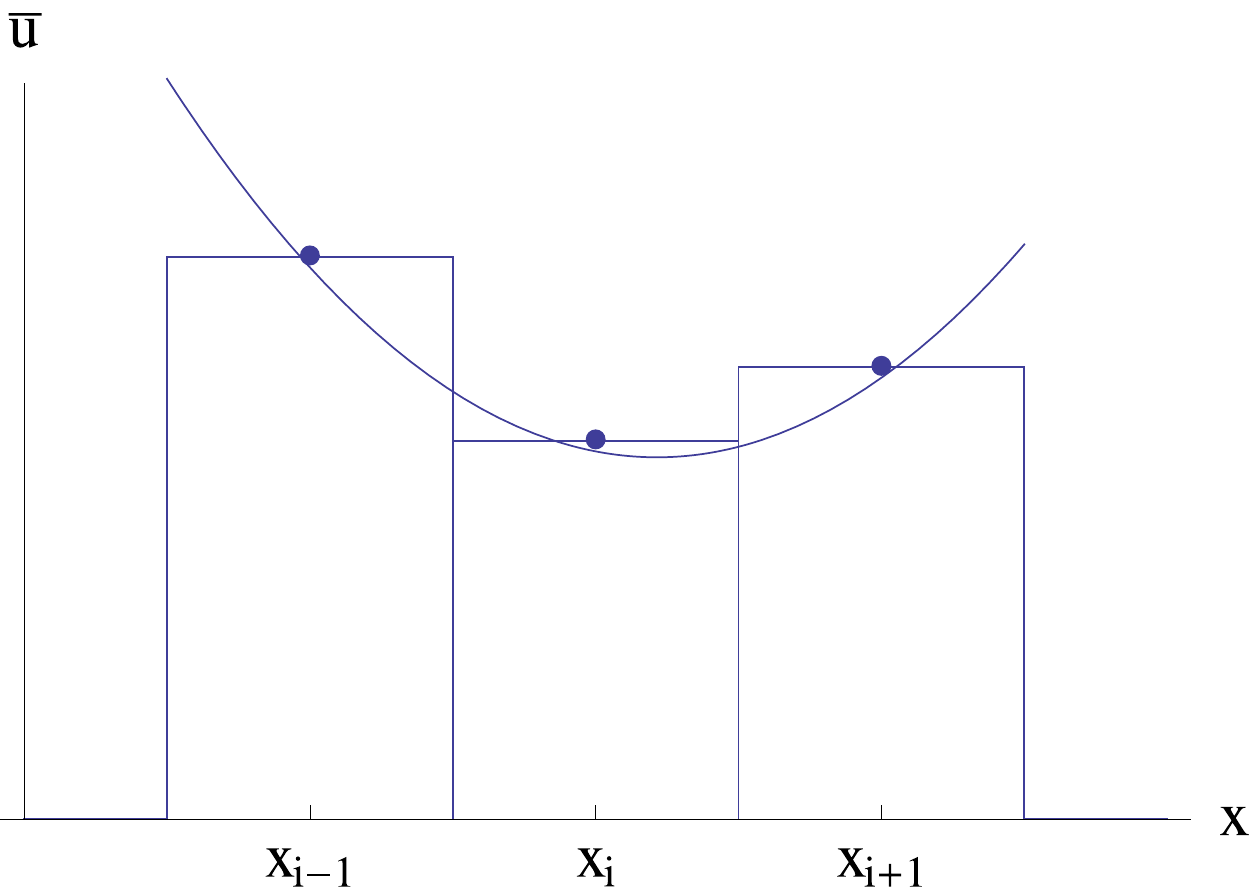}   		
	   		\caption{This situation is classified as smooth: $\tilde \phi_{\text{Lim}\mathcal{O}3} =\tilde \phi_{\mathcal{O}3}$.}
   		\end{subfigure}
   		\caption{Two similar situations that are treated differently by $\tilde \phi_{\text{Lim}\mathcal{O}3}$.}
   		\label{fig:unsymmetricDeltas}
   \end{figure}
	Note that the limiter function presented in \cite{CT2009} is not symmetric with respect to the diagonals. This means that for some cases 
	$\tilde \phi_{\text{Lim}\mathcal{O}3}(\delta_1,\delta_2) = \tilde \phi_{\mathcal{O}3}(\delta_1,\delta_2)$ but 
	$\tilde \phi_{\text{Lim}\mathcal{O}3}(-\delta_2, -\delta_1) \neq \tilde \phi_{\mathcal{O}3}(-\delta_2, -\delta_1)$, cf. Fig. \ref{fig:unsymmetricDeltas}. 
	This should not be the case. We therefore corrected this feature and defined the resulting limiter function $\tilde \phi_{\text{new}}$,
	\begin{align*}
		\tilde \phi_{\text{new}}(\delta_{i-1/2}, \delta_{i+1/2}) &= \phi_{\text{new}}(\theta_i)\;\delta_{i+1/2}, \\
		\phi_{\text{new}}(\theta_i) &= \max\left(0,\min\left( \tilde \phi_{\mathcal{O}3}, \max\left(-\theta_i,\min\left(2 \,\theta_i, \tilde \phi_{\mathcal{O}3}, 1.5 \right)\right)\right)\right).
	\end{align*}
	This new limiter function treats symmetric situations in the same manner, i.e. if $\tilde \phi_{\text{new}}(\delta_1,\delta_2) = \tilde \phi_{\mathcal{O}3}(\delta_1,\delta_2)$
	then also $\tilde \phi_{\text{new}}(-\delta_2, -\delta_1) = \tilde \phi_{\mathcal{O}3}(-\delta_2, -\delta_1)$.
\subsection{Decision Criterion}\label{subsec:theory2}
	On a three-point stencil, it is almost impossible to define a criterion that fully ascertains whether the function exhibits the beginning of a discontinuity or a smooth extremum. 
	As stated in Sec. \ref{subsec:theory1}, the two-parameter setting is the necessary prerequisite for the definition of such a criterion. In an earlier work \cite{CT2009}, {\v{C}}ada 
	and Torrilhon proposed a switch function $\eta(\delta_{i-1/2}, \delta_{i+1/2})$ which tests for smooth extrema. Their switch function defines an asymptotic region of radius $r$ 
	around the origin in the $\delta_{i+1/2} - \delta_{i-1/2} $ - plane in which we can safely switch to the third-order scheme. The limiter function $\phi_{\text{Lim}\mathcal{O}3}$
	together with this switch function has been successfully applied (e.g. \cite{K2011, KP2014, MTB2010, PKK2014}). Unfortunately, the authors do not specify the parameter $r$, which determines
	the size of the asymptotic region. With this idea in mind, we found that the most promising potential to distinguish discontinuities from smooth extrema is by measuring the 
	magnitude of the vector $(\delta_{i-1/2}, \delta_{i+1/2})$. When this vector is bounded in some appropriate norm, the reconstruction is switched to the full-third-order reconstruction, 
	even though one of lateral derivatives may be vanishing. 
	\begin{lemma}\label{lemma:magnitudeEta}
		In the vicinity of an extremum $\xi_0$, for $|x_i-\xi_0| \leq \Delta x$, the following relations hold:
		\begin{subequations}		
		\begin{align}
		\label{eq:lemma1}
			\begin{Vmatrix}			
				\begin{pmatrix}
					\delta_{i-1/2}\\ \delta_{i+1/2}
				\end{pmatrix}
			\end{Vmatrix}_2
			&\leq \,c\,\max_i |u_{0i}^{\prime\prime}|\,\Delta x^2\quad \text{with}\;c=\sqrt{\frac{5}{2}} 
			\\
			\label{eq:lemma2}
			\begin{Vmatrix}			
				\begin{pmatrix}
					\delta_{i-1/2}\\ \delta_{i+1/2}
				\end{pmatrix}
			\end{Vmatrix}_1
			&\leq \,c\,\max_i |u_{0i}^{\prime\prime}|\,\Delta x^2\quad \text{with}\;c=2
		\end{align}			
		\end{subequations}	
	\end{lemma}
	%
	%
	Lemma \ref{lemma:magnitudeEta} makes a statement on the magnitude of the differences across the cell interfaces. 
	The bound only depends on the grid size $\Delta x$ and the initial condition $u_0$. 
	\begin{definition}\label{def:eta}
		The switch function $\eta$ that marks the limit between smooth extrema and discontinuities is defined by
		\begin{align}
		\label{eq:eta}
 			\eta = \frac{\sqrt{\delta^2_{i-1/2}+\delta^2_{i+1/2}}}{\sqrt{\frac{5}{2}}\,\alpha\,\Delta x^2} \lessgtr 1	
 		\end{align}
 		with 
 		\begin{align}
 		\label{eq:alpha}
 			\alpha \equiv \max_{i \in \Omega \backslash \Omega_d} |u_{0_i}^{\prime\prime}(x)|.
 		\end{align} 		
 		Here, $\Omega$ is the computational domain and $\Omega_d$ is a set of points where the initial condition $u_0$ is discontinuous.
	\end{definition}
	\begin{proof} (Proof of Lemma \ref{lemma:magnitudeEta})\\
		Eq. \eqref{eq:lemma1} can be proven using a similar formulation of Def. \ref{def:eta}:
		\begin{align}
  	  		\frac{\delta^2_{i-1/2}+\delta^2_{i+1/2}}{(\alpha\,\Delta x^2)^2}  &= \frac{1}{\alpha^2}\left(\frac{u_{i+1} -2u_i+ u_{i-1}}{\Delta x^2}\right)^2 
		    +\frac{2}{\alpha^2 \Delta x^2}\left(\frac{u_{i+1}-u_i}{\Delta x}\right)\left(\frac{u_i-u_{i-1}}{\Delta x}\right)
		\end{align}
		A Taylor development around $x_i$ yields
		\begin{align}
  		\label{eq:etaBdiscMVlong}
    	\frac{\delta^2_{i-1/2}+\delta^2_{i+1/2}}{(\alpha\,\Delta x^2)^2}  &= \frac{1}{2}\left(\frac{u_i^{\prime\prime}}{\alpha}\right)^2
	    +\frac{2}{\Delta x^2}\left(\frac{u_i^\prime}{\alpha}\right)^2
	    + \frac{5}{6}\frac{u_i^\prime u_i^{(3)}}{\alpha^2}
   	    +\mathcal{O}(\Delta x^2).
  	 \end{align}
     In the vicinity of an extremum $\xi_0$, for $|x_i-\xi_0| \leq \Delta x$, the derivative fulfills $u_i^\prime \leq  u_{\xi_0}^{\prime\prime}\Delta x + \mathcal{O}(\Delta x^2)$. 
     Therefore, Eq. \eqref{eq:etaBdiscMVlong} reduces to
	\begin{align}
  		\label{eq:etaBdiscMV}	
		  \frac{\delta^2_{i-1/2}+\delta^2_{i+1/2}}{(\alpha\,\Delta x^2)^2} \leq \frac{1}{2}\left(\frac{u_i^{\prime\prime}}{\alpha}\right)^2
	    +2\left(\frac{u_{\xi_0}^{\prime\prime}}{\alpha}\right)^2+\mathcal{O}(\Delta x).
	\end{align}
    Setting $\alpha \equiv \max_{i \in \Omega \backslash \Omega_d} |u_{0_i}^{\prime\prime}(x)|$
	\begin{align*}
		\frac{\delta^2_{i-1/2}+\delta^2_{i+1/2}}{(\alpha\,\Delta x^2)^2} \leq \frac{5}{2}
	\end{align*}	     
     holds true, which shows Eq. \eqref{eq:lemma1}. \\
     In a similar manner, Eq. \eqref{eq:lemma2} can be proven.
	\end{proof}
	\noindent
	With Def. \ref{def:eta}, Lemma \ref{lemma:magnitudeEta} states that in the vicinity of smooth extrema, $\eta \leq 1$ holds. Combining this 
	information with the new limiter function $\tilde \phi_\text{new}$, we use this result to define the combined limiter
	%
	\begin{align*}
	 		&\tilde\phi_{\text{comb}}(\delta_{i-1/2},\delta_{i+1/2}) := 
	  		\begin{cases}
		  		\tilde\phi_{\mathcal{O}3}(\delta_{i-1/2},\delta_{i+1/2}) \quad\qquad &\text{if}\;\eta <1-\varepsilon\\
	  	  		\tilde\phi_{\text{new}}(\delta_{i-1/2},\delta_{i+1/2})\quad\; &\text{if}\;\eta > 1+\varepsilon\\
	  	  		W\left(\tilde\phi_{\mathcal{O}3}, \tilde\phi_{\text{new}}\right)\; &\text{else}.
	  		\end{cases}
	  \end{align*}   
	  where $\varepsilon$ is a small number of order $10^{-6}$ and $W(\cdot,\cdot)$ a linear function to ensure Lipschitz continuity 
	  of $\tilde\phi_{\text{comb}}$, cf. \cite{CT2009} for more details.
	  
\section{Numerical Results}\label{sec:numericalResults}
	In this section we want to test the decision criterion $\eta$ for the one-dimensional linear advection equation 
	\begin{subequations}	
		\label{eq:advectionEq}
		\begin{align}
			u_t + u_x &= 0, \\
			u(x,t=0) &= u_0(x), 
		\end{align}	
	\end{subequations}
	with two different characteristic initial conditions (ICs) on a periodic domain $[-1, 1]$. 
	Since $\eta$ requires the input of $\alpha = \max_{i\not \in \Omega_d} |u_{0_i}^{\prime\prime}|$ this external input is a possible 
	source of error. For this reason, we test for input values that are 
	\begin{enumerate}
		\item of the right order of magnitude
		\item over estimated, i.e too large
		\item under estimated, i.e. too small
		\item much too small.
	\end{enumerate}	
	The aim is to study the impact of possibly-incorrect input values and thus wrong switching functions $\eta$.
	\begin{figure}[ht!]
   	\centering
   		\begin{subfigure}{0.49\textwidth}
   			\centering
   			\includegraphics[width=\textwidth]{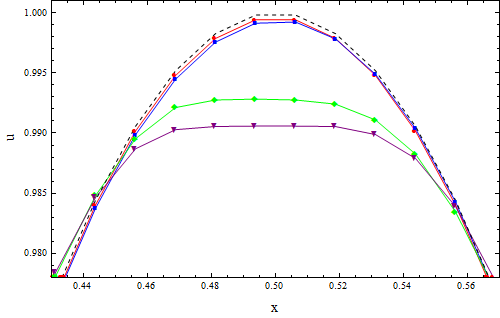}
			\caption{Solution of Eq. \eqref{eq:advectionEq} for different values of $\alpha$.}
			\label{fig:smoothDataSol}
   		\end{subfigure}
   		\hfill
   		\begin{subfigure}{0.49\textwidth}
   			\centering
   			\includegraphics[width=\textwidth]{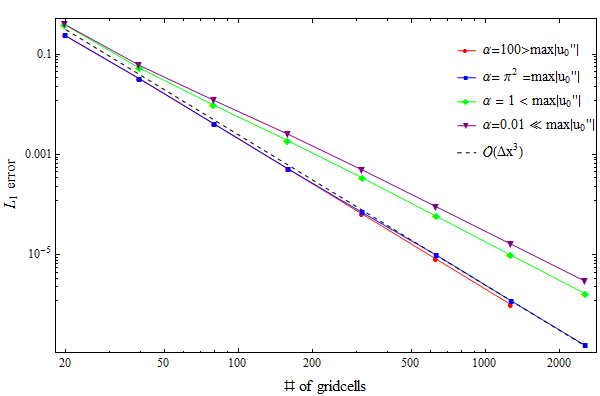}
	   		\caption{Double-logarithmic plot of the $\text{L}_1$-error vs. number of grid cells.}
			\label{fig:smoothDataError}
   		\end{subfigure}
		\caption{Results are calculated at $t_{\text{end}}=20.0$ using $u_0(x)=sin(\pi x), \nu=0.8, \Delta x = 0.0125$.}
		\label{fig:smoothData}
   \end{figure}
   %
\subsection{Convergence Studies For Smooth Initial Data}
   \noindent
	We solve the advection equation \eqref{eq:advectionEq} with the IC $u_0(x)=sin(\pi x),\, x\in [-1, 1]$. The function is convected until $\text{t}_\text{end}=20$ with 
	Courant number $\nu = 0.8$. In Fig. \ref{fig:smoothDataSol} we have plotted an area of interest of the solution of Eq. \eqref{eq:advectionEq}.  Fig. \ref{fig:smoothDataError} 
	shows the double-logarithmic $\text{L}_1$-error vs. number of grid cells. Both plots are depicted for different values of $\alpha$ and have been calculated for
	$n=160$ grid cells, i.e. $\Delta x = 0.0125$. \\
   	Fig. \ref{fig:smoothData} clearly points out that for the smooth test case, an over estimation of $\alpha$ does not effect the $3^{\text{rd}}$-order convergence 
   	of the solution. This is due to the fact that a large $\alpha$ means essentially no limiting but a direct application of the full $3^{\text{rd}}$-order reconstruction. 
   	If the input value 	for $\alpha$ is smaller, 
   	the limiter function $\tilde\phi_{\text{new}}$ is applied more often. In this case, a higher resolution is needed to distinguish between the discretization of a smooth 
   	extremum and a shallow gradient.
  	%
\subsection{Initial Condition with Discontinuous Data}
	%
	 In case of the square wave $u_0(x)=\mathds{1}_{[-0.5, 0.5]}(x),\, x \in [-1, 1]$, the input for $\alpha$, as defined by Eq. \eqref{eq:alpha} would yield $0$. 
	 However, this means that the new limiter function $\tilde \phi_{\text{new}}$ always takes effect and yields $0$ in most parts of the domain. This is because 
	 at least one of the consecutive slopes $\delta_{i\pm 1/2}=0$. However, arguing that in the smooth parts, $\delta_{i+1/2}\approx\delta_{i-1/2}$ (even though 
	 they yield $0$), we are close to the diagonals and thus, the $3^{\text{rd}}$-order reconstruction should be applied.
    \begin{figure}[ht]
		\begin{subfigure}{0.49\textwidth}
			\includegraphics[width=\textwidth]{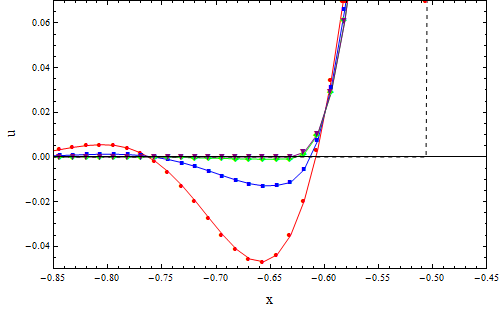}
			\caption{Solution of Eq. \eqref{eq:advectionEq} for different values of $\alpha$.}
		\end{subfigure}
		\begin{subfigure}{0.49\textwidth}
			\includegraphics[width=\textwidth]{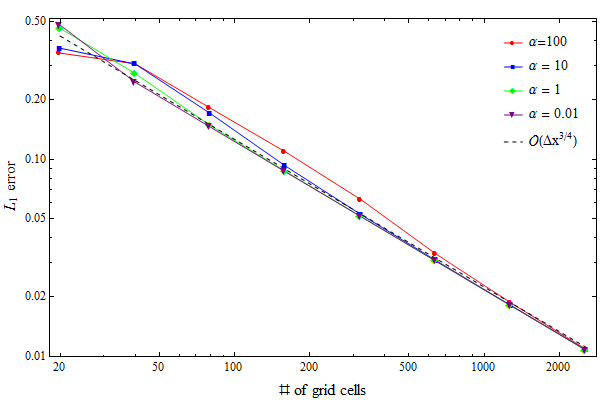}
			\caption{Double-logarithmic plot of the $\text{L}_1$-error vs. number of grid cells.}
		\end{subfigure}
		\caption{Results are calculated at $t_{\text{end}}=20.0$ using $u_0(x)=\mathds{1}_{[-0.5, 0.5]}(x)$, $\nu=0.8, \Delta x = 0.0125$.}
		\label{fig:discData}
	\end{figure}	
	Testing different values of $\alpha$ revealed that for larger values, the oscillatory behavior increases. This is due to the
	fact that with increasing $\alpha$ the region where $\tilde\phi_{\mathcal{O}3}$ is applied increases. Utilizing solely the 
	full $3^{\text{rd}}$-order reconstruction $\tilde\phi_{\mathcal{O}3}$ on the square wave is known to result in large over- and 
	undershoots and to asymptotically 	yield order $\mathcal{O}(\Delta x^{3/4})$.	
	Fig. \ref{fig:discData} shows that for small values of $\alpha$, the solution converges faster to $\mathcal{O}(\Delta x^{3/4})$ than 
	for large values of $\alpha$. Thus, small values should be preferred, however, even when the input for $\alpha$ is overestimated, 
	the solution converges when a sufficient number of grid cells is used.
	%
 		%

\end{document}